\documentclass[10pt]{article}      

\usepackage[dvips]{graphics}      
     
\usepackage{amssymb}      
\usepackage{amsmath}      
\usepackage{amsthm}      
\usepackage{amsbsy}      
\usepackage[dvips]{graphics}      
\usepackage{amscd}      
 
\theoremstyle{plain}      
\newtheorem{theorem}{Theorem}[section]

\newtheorem{lemma}{Lemma}[section]     
\newtheorem{corollary}{Corollary}[section]     
\newtheorem{proposition}{Proposition}[section]

\newtheorem{definition}{Definition}[section]      
\theoremstyle{remark}      
\newtheorem{remark}{Remark}[section]

\newcommand{\Z}{{\mathbb{Z}}}      
\newcommand{\C}{{\mathbb{C}}}      
\newcommand{\R}{{\mathbb{R}}}

\newcommand{\s}{{\cal S}}

\newcommand{\M}{{{\cal M}}}      
      
\newcommand{\B}{{{\cal B}}}

\newcommand{\T}{\widetilde{T}}

\newcommand{\su}{{\cal S}_{\infty}}      
\usepackage{amssymb}      
\def\N{\mathbb{N}}      
\def\C{\mathbb{C}}      
\def\Z{\mathbb{Z}}      
\def\Q{\mathbb{Q}}      
\def\T{\mathbb{T}}

\begin{document}

\title{An infinite genus mapping class group and stable cohomology\footnote{   This version: {\tt January 30, 2008}.  L.F. was partially supported by 
the ANR Repsurf:ANR-06-BLAN-0311.    
Available electronically at       
          \tt  http://www-fourier.ujf-grenoble.fr/\~{ }funar }}      
 \author{      
\begin{tabular}{cc}      
 Louis Funar &  Christophe Kapoudjian\\      
\small \em Institut Fourier BP 74, UMR 5582       
&\small \em Laboratoire Emile Picard, UMR 5580\\      
\small \em University of Grenoble I &\small \em University of Toulouse      
III\\      
\small \em 38402 Saint-Martin-d'H\`eres cedex, France      
&\small \em 31062 Toulouse cedex 4, France\\      
\small \em e-mail: {\tt funar@fourier.ujf-grenoble.fr}      
& \small \em e-mail: {\tt ckapoudj@picard.ups-tlse.fr} \\      
\end{tabular}      
}

\maketitle      
  
\begin{abstract} 
We exhibit a finitely generated group $\M$ whose rational homology is isomorphic to the 
rational stable homology of the mapping class group. It is defined as a mapping class
 group  associated to a surface $\su$ 
of infinite genus, and contains all the pure mapping 
class groups of compact surfaces of genus $g$ with $n$ boundary components, for any $g\geq 0$ and $n>0$.
We construct a representation of $\M$ into the 
restricted symplectic group ${\rm Sp_{res}}({\cal H}_r)$
of the real Hilbert space generated by the homology classes of non-separating 
circles on $\su$, which generalizes the classical symplectic representation of the mapping class groups. Moreover, we show 
that the first universal 
Chern class in $H^2(\M,\Z)$ is the pull-back of the Pressley-Segal 
class on the restricted linear group ${\rm GL_{res}}({\cal H})$ via the 
inclusion ${\rm Sp_{res}}({\cal H}_r)\subset {\rm GL_{res}}({\cal H})$. 
\vspace{0.2cm} 
 
\noindent 2000 MSC Classification: 57 N 05, 20 F 38, 22 E 65,81 R 10.  
 
\noindent Keywords: mapping class groups, infinite surface, Thompson group, 
stable cohomology,  Chern class, restricted symplectic group.  
\end{abstract}

\section{Introduction}  
\subsection{Statements of the main results}

The tower of all  extended mapping class groups  was considered first  by 
Moore and Seiberg (\cite{MS}) as part of the conformal field theory data.  
This object is actually a groupoid, which has been proved to be 
finitely presented (see \cite{BK,BK2,FG,HLS}).  When seeking for a group 
analog Penner (\cite{pe0})  investigated a  
universal  mapping class group which arises by means of a
completion process and which is closely related to the group of 
homeomorphisms of the circle, but it seems to be infinitely generated. \\

\noindent In \cite{FK}, 
we introduced  the universal mapping class group in genus zero ${\mathcal 
  B}$.  
 The latter is an extension of the  
Thompson's group  $V$ (see \cite{CFP}) by the infinite spherical pure 
mapping class group. We proved in \cite{FK} that the group $\B$ is finitely 
presented and we exhibited an explicit presentation.  Our main difference with  the previous attempts is that we consider 
groups acting on infinite surfaces with  a prescribed behaviour at infinity
that comes from actions on trees.  \\

\noindent Following the same kind of approach, we propose a treatment of the arbitrary  genus case by introducing
 a mapping class 
group $\M$, called {\it the asymptotic infinite genus mapping class group}, that contains  a large part of the mapping class groups 
of compact surfaces with boundary. More precisely, the group $\M$  contains all 
the pure mapping class groups $P\M(\Sigma_{g,n})$ of compact surfaces    
$\Sigma_{g,n}$ of genus $g$ with $n$ boundary components, for any $g\geq 0$ and $n>0$.  Its construction is roughly as follows.
Let ${\cal S}$ denote the surface 
obtained by taking the boundary of the  3-dimensional 
thickening of the complete trivalent tree, and  
further let ${\cal  S}_{\infty}$ be the result of attaching a handle 
to each cylinder in ${\cal S}$ that corresponds to an edge of the tree 
(see figure 1).  Then $\M$ is the group of mapping classes of those homeomorphisms 
of  ${\cal  S}_{\infty}$ which preserve a certain {\it rigid structure} at infinity (see Definition \ref{defm} for the precise definition). 
This rigidity condition essentially implies that $\M$ induces a group of transformations on the set
of ends of the tree, which is isomorphic to Thompson's group $V$. The relation between both groups is enlightened 
by a short exact sequence $1\rightarrow P\M\rightarrow \M\rightarrow V\rightarrow 1$,
 where $P\M$ is the mapping class group of compactly supported homeomorphisms of $\su$. The latter is
 an infinitely generated group. Our first result is:     
 \begin{theorem}\label{fg}     
The group $\M$ is finitely generated.      
\end{theorem}   

\noindent The interest in considering the group $\M$, outside the      
framework of the topological quantum field theory where it can replace     
the duality groupoid, is the following homological property:

\begin{theorem}\label{coho}      
The rational homology of $\M$ is isomorphic to the stable rational     
homology of the (pure) mapping class groups.      
\end{theorem} 

\noindent As a corollary  of the argument of the proof (see Proposition \ref{c1}), the group $\M$ is perfect, and 
$H_2(\M,\Z)=\Z$. For a reason that will become clear in what follows, the generator of $H^2(\M,\Z)\cong \Z$ is called 
{\em the first universal Chern class} of $\M$, and is denoted $c_1(\M)$. \\

\noindent Let ${\cal M}_g$ be the mapping class group of a closed surface $\Sigma_g$ of genus $g$. We show that  the standard representation 
$\rho_g: {\cal M}_g\rightarrow {\rm Sp}(2g,\Z)$ in the symplectic group,  
deduced from the 
action of ${\cal M}_g$ on  $H_1(\Sigma_g,\Z)$, extends to 
the infinite genus case, by replacing the finite dimensional setting by concepts of Hilbertian analysis. In particular, a key role is played by Shale's {\it restricted 
  symplectic group} ${\rm Sp_{res}}({\cal H}_r)$ on the real Hilbert space ${\cal 
  H}_r$ generated by the homology classes of non-separating closed curves of ${\cal 
  S}_{\infty}$. We have then:  
 
\begin{theorem}\label{metab}
The action of $\M$ on $H_1({\cal S}_{\infty},\Z)$ induces a representation $\rho:\M\rightarrow {\rm Sp_{res}}({\cal H}_r)$. 
\end{theorem}  

\vspace{0.2cm} \noindent
 The generator $c_1$ of $H^2({\cal M}_g,\Z)$ is called the 
first Chern class, since it may be obtained as follows 
(see, e.g., \cite{mo}). The group ${\rm Sp}(2g,\Z)$ is 
contained in the symplectic group ${\rm Sp}(2g,\R)$, whose maximal compact subgroup 
is the unitary group $U(g)$. Thus, the first Chern class may be viewed in 
$H^2(B {\rm Sp}(2g,\R),\Z)$. It can be  first pulled-back on 
$H^2(B{\rm Sp}(2g,\R)^{\delta},\Z)=H^2({\rm Sp}(2g,\R),\Z)$ and 
then on $H^2({\cal   M}_g,\Z)$ via $\rho_g$. 
This is the generator of $H^2({\cal M}_g,\Z)$. 
Here $B{\rm Sp}(2g,\R)^{\delta}$ 
denotes the classifying space of the group ${\rm Sp}(2g,\R)$ endowed with  
the discrete topology. \\

 \noindent The restricted symplectic group ${\rm Sp_{res}}({\cal H}_r)$ has a well-known 
2-cocycle, which measures the 
projectivity of the {\it Berezin-Segal-Shale-Weil metaplectic representation} 
in the bosonic Fock space (see \cite{ne}, Chapter 6 and Notes 
p. 171). Contrary to the finite dimension case, this cocycle is not directly related 
to the topology of ${\rm Sp_{res}}({\cal H}_r)$, 
since the latter is a contractible Banach-Lie group. 
However, ${\rm Sp_{res}}({\cal H}_r)$ embeds into the 
restricted linear group of Pressley-Segal ${\rm GL^0_{res}}({\cal H})$ (see 
\cite{pr-se}), where ${\cal 
  H}$ is the complexification of ${\cal H}_r$, which possesses a cohomology class of degree 2: the Pressley-Segal 
class $PS\in H^2({\rm GL^0_{res}}({\cal H}),\C^*)$. The group ${\rm GL^0_{res}}({\cal H})$ 
is a homotopic model of the classifying space $BU$, where 
$U=\displaystyle{\lim_{n\to \infty} U(n,\C)}$, and the class $PS$ does correspond to the 
universal first Chern class. Its restriction on ${\rm Sp_{res}}({\cal H}_r)$ is 
closely related to the Berezin-Segal-Shale-Weil cocycle, and reveals the topological 
origin of the latter. Via the composition of morphisms 
$$\M\longrightarrow {\rm Sp_{res}}({\cal H}_r) \hookrightarrow 
{\rm GL^0_{res}}({\cal 
  H}),$$  we then derive from 
$PS$ an {\it integral} cohomology class on $\M$ (see Theorem \ref{chernimp}
for a more precise statement): 
 
\begin{theorem}\label{chern} 
The Pressley-Segal class $PS\in H^2({\rm GL^0_{res}}({\cal H}),\C^*)$ induces the first universal Chern 
class $c_1(\M)\in H^2({\cal M},\Z)$. 
\end{theorem}

\noindent      
{\bf Acknowledgements.} The authors are indebted to Vlad Sergiescu 
for enlighting discussions and particularly for suggesting the existence 
of a connection between the first universal Chern class of ${\cal M}$ 
and the Pressley-Segal class. They are thankful to the referees 
for suggestions improving the exposition.

\subsection{Definitions}

\subsubsection{The infinite genus mapping class group ${\cal M}$}     
Set ${\M}(\Sigma_{g,n})$  for the extended 
mapping class group of the       
$n$-holed orientable surface  $\Sigma_{g,n}$ of genus $g$,      
consisting of the isotopy classes of orientation-preserving 
homeomorphisms      
of $\Sigma_{g,n}$ which respect a fixed parametrization of the 
boundary circles, allowing them to be permuted among themselves.

\vspace{0.2cm}
\noindent 
We wish to construct a mapping class group, containing      
all mapping class groups ${\M}(\Sigma_{g,n})$. 
It seems impossible to construct such a group, 
but if one relaxes  slightly our requirements then we could follow 
our previous method used for the genus zero case in \cite{FK}. 

\vspace{0.2cm}
\noindent 
The choice of the extra structure involved in the 
definitions below is important because the final result might depend on it.  
For instance, using the same planar punctured surface but different  
decompositions one obtained in \cite{FK2} 
two non-isomorphic braided Ptolemy-Thompson groups.

\begin{definition}[The infinite genus surface $\s_{\infty}$]\label{surfinfgenus}
Let ${\mathcal T}$ be the complete trivalent planar tree and 
$\s$ be the surface obtained by taking the boundary of the  
3-dimensional thickening of  ${\mathcal T}$.

\vspace{0.2cm}
\noindent 
By grafting
an {\em edge-loop} (i.e. the graph obtained by attaching a loop to a boundary vertex of an edge) 
at the midpoint of each edge of ${\cal T}$, one obtains the graph ${\cal T}_{\infty}$. The surface $\s_{\infty}$ is the boundary of the 3-dimensional thickening of ${\cal T}_{\infty}$.     
\end{definition}

\noindent 
The graph ${\mathcal T}$ (respectively ${\cal T}_{\infty}$) is 
embedded in $\s$ (respectively  $\s_{\infty}$) as a cross-section of 
the fiber projection, as indicated on figure 1. 
Thus,  $\su$ is obtained by removing 
small disks  from $\s$ centered at midpoints of edges 
of ${\mathcal T}$ and gluing back one holed 
tori $\Sigma_{1,1}$, called {\em wrists} which correspond to the 
thickening of edge-loops.

\noindent  It is convenient to assume
that 
${\cal T}$ is embedded in a horizontal plane, while the edge-loops are in vertical planes (see figure 1).

%

\begin{definition}[Pants decomposition of $\s_{\infty}$]        
         
A {\em pants decomposition} of the surface $\su$ is a maximal 
collection       of distinct nontrivial simple closed curves 
on $\su$ which are pairwise            disjoint and non-isotopic. 
The complementary regions (which are 3-holed           
spheres) are called {\em pairs of pants}.           

\noindent 
By construction, $\su$ is naturally equipped with a pants           
decomposition, which will be referred to below as the            
{\em canonical (pants) decomposition}, as shown in figure 1: 
\begin{itemize}       
\item the wrists  are decomposed using a meridian circle and the boundary  
circle of $\Sigma_{1,1}$. 
\item there is one pair of pants for each edge, which has one  
boundary circle for attaching the wrists, and two circles to grip to the 
other type of pants. We call them {\em edge pants}.   
\item there is one pair of pants for each vertex of the tree, called 
 {\em  vertex pants}.  
\end{itemize}

\begin{figure}
\begin{center}      
\includegraphics{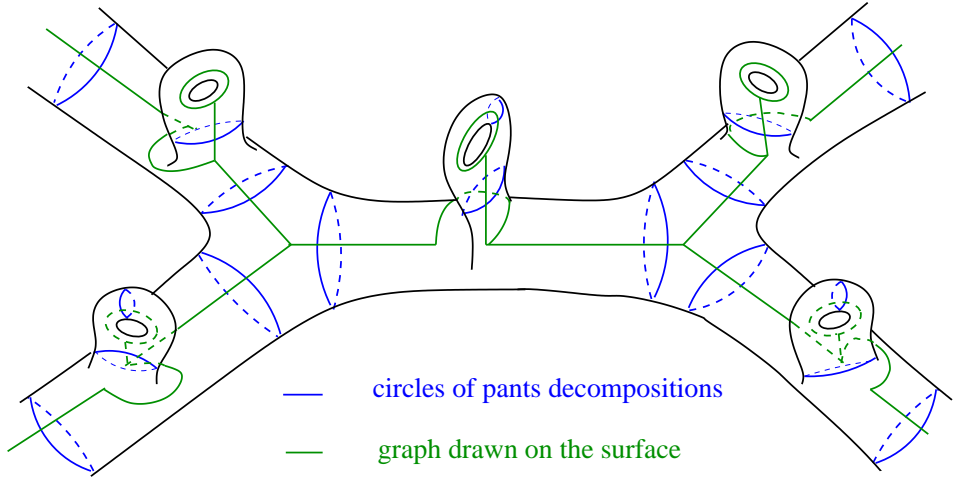}      
\caption{The infinite genus surface $\su$ with its canonical rigid structure}
\end{center}      
\end{figure}
\end{definition}

\noindent 
A pants decomposition is {\em asymptotically  trivial} if outside a compact subsurface of $\su$, it coincides with the           
  canonical pants decomposition.

\begin{definition}\label{defm}
\begin{enumerate}

\item A connected subsurface $\Sigma$ of $\su$ is admissible 
if all its boundary circles are 
from vertex type pair of pants from the canonical decomposition and 
moreover, if one boundary circle  from a vertex type pants 
is contained in $\Sigma$ then the entire pants is contained in $\Sigma$.  
In particular, $\su-\Sigma$ has no compact 
components.  

\item  Let $\varphi$ be a homeomorphism of $\su$. One says that       
$\varphi$ is {\em asymptotically rigid} if the following conditions are      
fulfilled:      
\begin{itemize}      
\item There exists an admissible subsurface $\Sigma_{g,n}\subset \su$ such     that $\varphi(\Sigma_{g,n})$ is also admissible.       
\item The complement $\su -\Sigma_{g,n}$ is a union of $n$ infinite      
  surfaces. Then the restriction $\varphi: \su -\Sigma_{g,n}\to \su-\varphi(\Sigma_{g,n})$      
 is {\em rigid}, meaning that it maps the pants      
 decomposition into the pants decomposition and 
maps  ${\cal T}_{\infty}\cap (\su -\Sigma_{g,n})$
 onto  ${\cal T}_{\infty}\cap (\su -\varphi(\Sigma_{g,n}))$. 
Such a surface $\Sigma_{g,n}$ is called a {\em support} for $\varphi$.      
\end{itemize}      
One denotes by $\M=\M(\su)$ the {\em group of asymptotically rigid      
  homeomorphisms of $\su$} up to isotopy  and call it the 
  {\em asymptotic mapping class group of infinite genus}.      
\end{enumerate}
\end{definition} 
\noindent In the same way one defined the asymptotic  mapping class group 
$\M(\s)$, denoted by $\B$ in \cite{FK}.

\begin{remark}
In genus zero (i.e. for the surface $\s$) a homeomorphism between two 
complements of admissible subsurfaces which maps the 
restrictions of the tree $\mathcal T$ one into the other is  
rigid, thus preserves the isotopy class of the pants decomposition. 
This is not anymore true in higher genus: the Dehn twist along a longitude 
preserves the edge-loop graph but it is not rigid, as a homeomorphism 
of the holed torus. 
\end{remark}

\begin{remark}
Notice that, in general,  rigid homeomorphisms $\varphi$  
do not have an invariant support i.e. an admissible 
$\Sigma_{g,n}$ such that $\varphi(\Sigma_{g,n})=\Sigma_{g,n}$. 
Take for instance a homeomorphism which translates the wrists 
along a geodesic ray in $\mathcal T$. 
\end{remark}

\begin{remark}      
Any admissible subsurface $\Sigma_{g,n}\subset \su$ has $n=g+3$.
Moreover $\su$ is the ascending union 
$\cup_{g=1}^{\infty}\Sigma_{g,g+3}$.        
Instead of the wrist $\Sigma_{1,1}$ use a surface  of higher genus $\Sigma_{g,1}$ and 
the same definitions as above. The admissible subsurfaces will be 
$\Sigma_{kg, k+3}$. 
The asymptotic mapping class group obtained this way is finitely generated 
by small changes in the proof below. 
\end{remark}

\begin{remark} 
The surface $\su$ contains infinitely many
compact surfaces of type $(g,n)$ with at least one boundary component. For any such compact subsurface      
$\Sigma_{g,n}\subset \su$, there is an obvious injective morphism       
$i_*:{P\M}(\Sigma_{g,n})\hookrightarrow P\M\subset \M$.  However, the morphism $i_*:{\M}(\Sigma_{g,n})\hookrightarrow \M$ 
is not always defined. Indeed, it exists if and only if the $n$ connected components of $\su \setminus \Sigma_{g,n}$
are homeomorphic to each other, by asymptotically rigid homeomorphisms.\\
In particular, for any admissible subsurface $\Sigma_{g,n}$ (hence $n=g+3$), $i_*$ extends to an injective morphism
$i_*:{\M}(\Sigma_{g,n})\hookrightarrow \M$ defined by rigid extension of homeomorphisms of $\Sigma_{g,n}$
to $\su$.
\end{remark}

\subsubsection{The group ${\mathcal M}$ and the Thompson groups}

\begin{definition}
\begin{enumerate}
\item Let ${\cal T}$ be the planar trivalent tree. A {\em partial 
tree automorphism} of ${\cal T}$ is an isomorphism of graphs $\varphi:{\cal T}\setminus \tau_{1}
\rightarrow {\cal T}\setminus \tau_{2}$, where $\tau_{1}$ and $\tau_{2}$ are two finite
trivalent subtrees of ${\cal T}$ (each vertex except the leaves are 3-valent). A connected component 
of ${\cal T}\setminus \tau_{1}$  or ${\cal T}\setminus \tau_{2}$
is a {\em branch}, that is, a rooted planar binary tree whose vertices are 3-valent, except the root, which is
2-valent. Each vertex of a branch has two descendant edges, and given an orientation to the plane,
one may distinguish between the left and the right descendant edges. A partial automorphism  $\varphi:{\cal T}\setminus \tau_{1}
\rightarrow {\cal T}\setminus \tau_{2}$ is {\em planar} if it maps each branch of  ${\cal T}\setminus \tau_{1}$  
onto the corresponding branch of ${\cal T}\setminus \tau_{2}$ by respecting the left and right ordering of the edges.
\item Two planar partial automorphisms $\varphi:{\cal T}\setminus \tau_{1}
\rightarrow {\cal T}\setminus \tau_{2}$ and $\varphi':{\cal T}\setminus \tau'_{1}
\rightarrow {\cal T}\setminus \tau'_{2}$ are {\it equivalent}, which is denoted $\varphi\sim \varphi'$, if and only if there exists a third 
$\varphi'':{\cal T}\setminus \tau''_{1}
\rightarrow {\cal T}\setminus \tau''_{2}$ such that $\tau_{1}\cup  \tau'_{1}\subset \tau''_{1}$,  $\tau_{2}\cup  \tau'_{2}\subset \tau''_{2}$ and $\varphi_{{\mid }{\cal T}\setminus \tau''_{1}}=
\varphi'_{{\mid }{\cal T}\setminus \tau''_{1}}=\varphi''$.
\item If $\varphi$ and $\varphi'$ are planar partial automorphisms, one can find $\varphi_{0}\sim \varphi$ and $\varphi'_{0}\sim \varphi'$  such that the source of $\varphi_{0}$ and the target of 
$\varphi'_{0}$ coincide. The product $[\varphi]\cdot [\varphi']=[\varphi_{0}\circ\varphi'_{0}]$ is well defined,
 as is easy to check. The set of equivalence classes of such automorphisms  endowed with
the above internal law, is a group with neutral element the class of $id_{\cal T}$. This is the {\em Thompson group $V$}.
\end{enumerate}
\end{definition}


\begin{remark}      
We warn the reader that our definition of the group $V$ is different      
from the standard one (as given in \cite{CFP}). Nevertheless, the present group $V$ is isomorphic to      
the group denoted by the same letter in \cite{CFP}.       
\end{remark}       
      
\noindent We introduce Thompson's group $T$, the subgroup of $V$ acting on the      
circle (see \cite{GS}), which will play a key role in the proofs.          
          
\begin{definition}[Ptolemy-Thompson's group $T$]\label{tho} 
Choose a vertex $v_{0}$ of $\cal T$. Each $g\in V$ may be represented by a planar partial automorphism  
 $\varphi:{\cal T}\setminus \tau_{1}\rightarrow {\cal T}\setminus \tau_{2}$   such that $v_{0}$
   belongs to $\tau_{1}\cap \tau_{2}$. Let $D_{1}$ (respectively $D_{2}$) be a disk containing 
   $\tau_{1}$ (respectively $\tau_{2}$), whose boundary circle $S_{1}$ (respectively $S_{2}$) passes through the leaves of $\tau_{1}$ 
   (respectively $\tau_{2}$), giving to them a cycling ordering. If $\varphi$ preserves this cycling ordering, which amounts to saying that the 
   bijection from the set of leaves of ${\tau}_{1}$ onto the set of leaves of ${\tau}_{2}$ can be extended to an orientation preserving homeomorphism
   from $S_{1}$ onto $S_{2}$, then any other  $\varphi'$ equivalent to  $\varphi$ also does, and one says that $g$ itself is {\em circular}. 
   The subset of circular elements of $V$ is a subgroup, called the {\em Ptolemy-Thompson group $T$}.    
\end{definition}

\begin{proposition}\label{sequence}
Set $P\M$ for the inductive limit of the pure mapping class groups      
of admissible subsurfaces of $\su$.       
We have then the following exact sequences:      
\[ 1 \to P\M \to \M \to V\to 1.\]              
\end{proposition}           
           
\begin{proof}  
         
Let $\varphi$ be an asymptotically rigid homeomorphism of $\su$ and $\Sigma_{g,n}$ a support 
for $\varphi$. Then it maps ${\cal T}_{\infty}\cap (\su -\Sigma_{g,n})$
 onto  ${\cal T}_{\infty}\cap (\su -\varphi(\Sigma_{g,n}))$, hence  
${\cal T}\cap (\su -\Sigma_{g,n})$
 onto  ${\cal T}\cap (\su -\varphi(\Sigma_{g,n}))$ by forgetting the action on the edge-loops. This may be identified with a planar partial
 automorphism $\phi:{\cal T}\setminus \tau_{1}\rightarrow {\cal T}\setminus \tau_{2}$. The map $[\varphi]\in \M\rightarrow [\phi]\in V$ is a group
 epimorphism. The kernel is the subgroup of           
isotopy classes of homeomorphisms inducing the identity outside a   
support,   
and hence is the direct limit of the pure mapping class groups.   
\end{proof}
      
\begin{remark} 
In \cite{FK} we prove the existence of a similar short exact sequence relating $\B$ to $V$, which splits over the Ptolemy-Thompson group $T$. It is worth noticing that the present extension of $V$ is not split 
over $T$. 
\end{remark}

\section{The proof of theorem \ref{fg}}

\subsection{Specific elements of  $\M$}

\noindent 
Recall that $\su$ has a canonical pants decomposition, as shown 
in figure 1. We fix an admissible subsurface $A=\Sigma_{1,4}$ which contains a 
central wrist and an admissible $B=\Sigma_{0,3}\subset \Sigma_{1,4}$ which is not 
adjacent to the wrist. 

\vspace{0.2cm} 
\noindent
Let us consider now the elements of $\M$ described in the pictures 
below. Specifically:       
\begin{itemize}      
\item Let $\gamma$ be a circle contained inside $B$ and parallel to the boundary curve labeled 3. Let $t$ be the right Dehn twist around $\gamma$. This means that, given an
outward orientation to the surface, $t$ maps an arc crossing $\gamma$ transversely to an arc   
  which turns right as it approaches $\gamma$. The dashed arcs (also called {\it seams}) on the left hand side picture figure out
  the boundary of the visible side of $B$. Their images by $t$ are represented on the right hand side picture,        
\begin{center}
\includegraphics{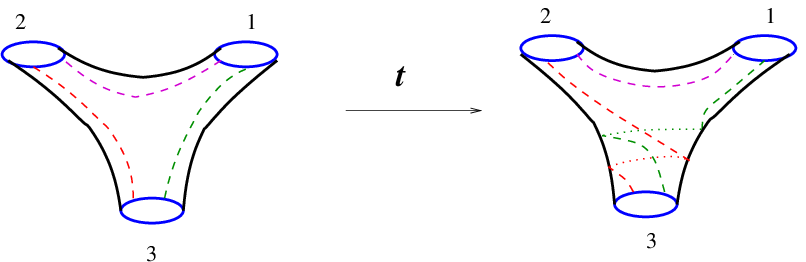}  
\end{center}
      
\item $\pi$ is the braiding, acting as a braid in      
  $\M(\Sigma_{0,3})$, with the support $B$. 
It rotates the circles 1 and 2 in the      
  horizontal plane (spanned by the circles) counterclockwise.

Assume that $B$  is identified with   
  the complex domain $\{|z|\leq 7, |z-3|\geq 1, |z+3|\geq 1\}\subset \C$.    
  A specific homeomorphism in the mapping class of $\pi$    
  is the composition of the counterclockwise   
 rotation of $180$ degrees around the origin --- which exchanges   
the small boundary circles labeled 1 and 2 in the figure ---  
with a map which rotates of $180$ degrees in the clockwise direction   
each boundary circle. The latter can be constructed as follows.     
  
\vspace{0.1cm}  
\noindent   
Let $A$ be an annulus in the plane, which we suppose for simplicity   
to be $A=\{1\leq |z|\leq 2\}$.    
The homeomorphism $D_{A, C}$ acts as the counterclockwise  
rotation of $180$ degrees  
on the boundary circle $C$ and   
keeps the other  boundary   
component pointwise   
fixed:   
\[D_{A, C}(z)=   
\left\{ \begin{array}{ll}  
z \exp(\pi\sqrt{-1}(2-|z|)), & \mbox{ if }  C=\{|z|=1\}\\  
z \exp(\pi\sqrt{-1}(|z|-1)), & \mbox{ otherwise}  
\end{array}\right.\]

\vspace{0.1cm}  
\noindent   
The map we wanted is    
$D^{-1}_{A_0, C_0}D^{-1}_{A_1, C_1}D^{-1}_{A_2, C_2}$, where   
$A_0=\{6\leq |z|\leq 7\}$, $C_0=\{|z|=7\}$,   
$A_1=\{1\leq |z-3|\leq 2\}$, $C_1=\{|z-3|=1\}$,   
$A_2=\{1\leq |z+3|\leq 2\}$, and  $C_2=\{|z+3|=1\}$.    
   
\noindent    
One has pictured also the images of the seams.  
\begin{center}
\includegraphics{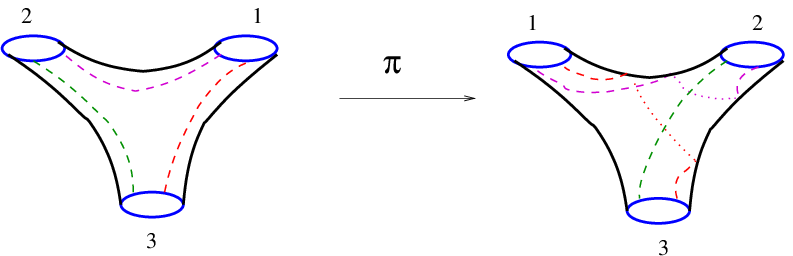}  
\end{center}
    
\item $\beta$ is the order 3 rotation in the vertical plane of the paper.
 It is the unique globally rigid   
  mapping class which permutes counterclockwise and cyclically the three   
  boundary circles of $B$. An invariant support for $\beta$ is $B$. 

\begin{center}
\includegraphics{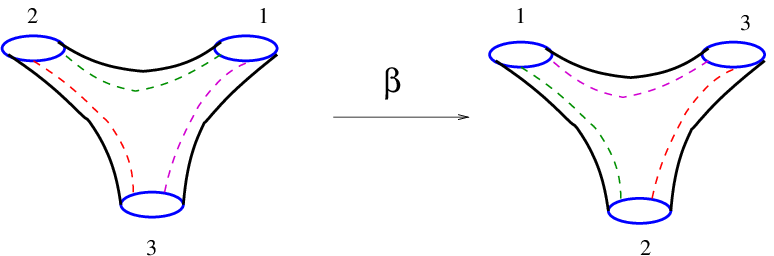}  
\end{center}
      
\item $\alpha$ is a twisted rotation of order 4 in the vertical 
plane which moves cyclically the labels of the boundary circles      
counterclockwise. Its support is a 4-holed torus $A=\Sigma_{1,4}$. 

\begin{center}
\includegraphics{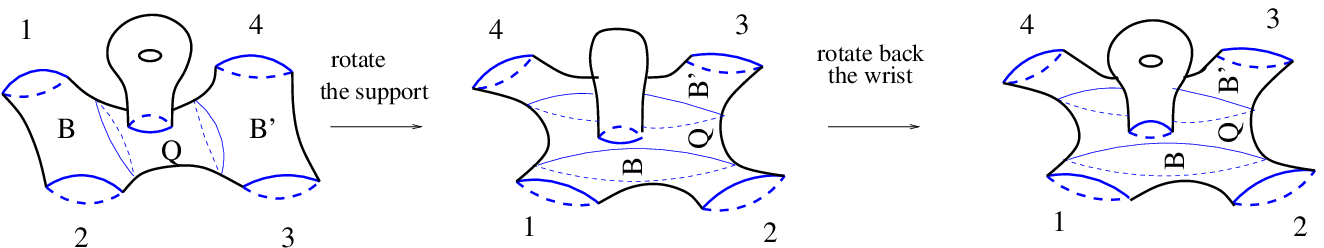}  
\end{center}
Let $\Sigma_{0,5}$ be the 5-holed sphere    
consisting of the union of  $B$  with the edge pants $Q$ near $B$ 
and the next vertex pants $B'$ adjacent to $Q$.
There are four boundary circles which are vertex type and one boundary 
circle which bounds a wrist. We perform first a rotation in $\R^3$  
which   preserves globally the pants decomposition and visible side, 
permutes counterclockwise and   
  cyclically the four vertex type boundary circles  
of $\Sigma_{0,5}$ and rotates the edge type circle according to 
one fourth twist. This rotation changes the  position of the 
wrist $\Sigma_{1,1}$ in $\R^3$. 
We consider next the clockwise rotation of this wrist alone, of  
angle $\frac{\pi}{2}$ around the vertical axis  that meets 
the edge type circle in its center.  
This rotation restores the initial wrist position. 
The composition of the two partial rotations above is a homeomorphism 
of $\Sigma_{1,4}$ that gives a well-defined element of $\M$.

\item  Let $a_1, b_1$ be the  meridian and longitude on the basic wrist 
in $A$. We denote by $t_{a_1}$ and $t_{b_1}$ 
the Dehn twists along these curves. Further, 
$t_0$ states for the  Dehn twist $t_0$ along the boundary 
circle of the wrist. 
\end{itemize}      
      
\begin{remark}         
It is worthy to note that we have       
three types of Dehn twists: those along separating curves       
(conjugate either to $t$ (the boundary on the vertex pants) or with 
$t_0$ (the edge type pants) and those along non-separating curves which  
are conjugate to  the twist around such a curve on the wrist.      
\end{remark}

\subsection{Generators for  $P\M$} 

\noindent Consider the following collection of simple curves drawn on 
$\su$:

\begin{center}
\includegraphics{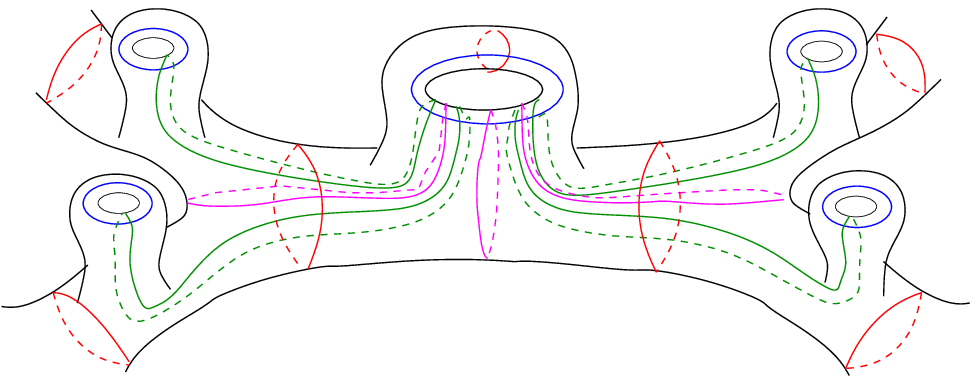}
\end{center}

\noindent 
Their description follows. 
\begin{enumerate}
\item Choose, for each wrist, a longitude $b_i$, 
which turns once along  the wrist. 
\item For each pair of wrists  we choose a circle
joining them as follows. For each wrist we have an  
arc going  from the base point of its attaching circle to the longitude 
and back to the opposite point of the circle.  
Then join these two pairs of points by a pair of 
parallel arcs in the horizontal surface, asking that the arc which joins the two base points be a geodesic path
in the tree ${\cal T}$. 
We call them {\em wrist-connecting loops}. 
\item Further we associate a loop to each pair consisting of a wrist and a vertex of the tree $\cal T$. 
A vertex gives rise to a pair of pants in $\su $. Two of the boundary components of these pants 
correspond to the directions to move away from the wrist. 
Thus we can define again an arc on the pair of pants 
which joins a point $p$ of the third circle (closest to the wrist), on the visible side of $\su$,
to its opposite, on the hidden side, and separates the remaining two circles. 
Consider the loop resulting from gathering the following three kinds of arcs: 
\begin{enumerate}
\item the arc on the  
wrist;
\item the arc on the pair of pants;
\item  and a pair of parallel arcs 
which join them, asking that the arc which joins the point $p$ to the base point of the wrist be a geodesic path
of the tree ${\cal T}$.
\end{enumerate}
 We call them 
{\em vertex connecting loops}. 
\item Consider the loops that come from the  
canonical pants decomposition of $\s$ by doubling them. 
We call them the {\em horizontal pants decomposition loops}.
\end{enumerate}


\begin{lemma}
The set $\cal H$ of Dehn twists along the meridians, the longitudes, the wrist 
connecting loops  associated to edges, the 
vertex connecting loops and the horizontal pants decomposition loops
generates $P\M$.  
\end{lemma}
\begin{proof}
It suffices to consider the  finite case of an 
admissible surface with boundary that contains $g$ 
wrists and has $g+3$ boundary components. 
Then the lemma follows from \cite{Ge}, in which it is proved that the pure mapping class group
of such a surface is generated by a set of Dehn twists ${\cal H}_{g,n}$ (with $n=g+3$). It suffices
to check that all the Dehn twists belonging to  ${\cal H}_{g,n}$ also belong to the set $\cal H$.
Referring to the notations of \cite{Ge}, there are four types of Dehn twists in ${\cal H}_{g,n}$: the $\alpha_{i}$'s, the
$\beta_{i}$'s, the $\gamma_{ij}$'s and the $\delta_{i}'s$. The $\alpha_{i}$'s are associated to 
wrist-connecting or vertex connecting loops, the $\beta_{i}$'s are associated to longitudes $b_{i}$'s, the $\gamma_{ij}$'s 
are associated to wrist-connecting loops, except $\gamma_{12}$ which is associated to a vertex-connecting
loop, and finally, the $\delta_{i}$'s are associated to the circles of the boundary of the surface, hence of the pants decomposition 
(after doubling them) of $\su$.
Therefore all of them belong to $\cal H$.

\vspace{0.2cm}\noindent 
Remark also that in (\cite{Ge}, figure 1) the 1-handles are 
cyclically ordered  and arranged on one side 
and  then followed by all boundary components 
of the surface. However, we can  arbitrarily permute 
the position of holes and 1-handles  in the picture and 
keep the same system of generators.  
\end{proof} 
 
\subsection{The action of $\T$ on the generators of $P\M$}

\subsubsection{The groups $\T$ and $T^*$}

Consider the subgroup $\T$ of $\M$ generated by the elements $\alpha$ and 
$\beta$. We will prove that the set of conjugacy classes for the action of $\T$ on
$\cal H$ is finite by considering the action of $\T$ on some planar 
subsurface of $\su$.  

\vspace{0.2cm}\noindent 
The surface obtained by puncturing (respectively deleting disjoint 
small open disks from) $\s$ at the midpoints  of the edges 
is denoted by $\s^*$ (and respectively $\s^{\bullet}$). 
The 2-dimensional thickening in $\s$  
of the embedded tree ${\mathcal T}$ 
is an infinite planar surface, which will be called the {\em visible side} of $\s$, and will be denoted $D$. The intersection of $D$ with 
$\s^*$ and $\s^{\bullet}$ is denoted $D^*$ and $D^{\bullet}$, respectively. 

\vspace{0.2cm}\noindent 
The elements $\alpha$ and $\beta$ as defined above (i.e. as 
specific homeomorphisms, not only as mapping classes)  keep 
invariant both $\s^{\bullet}$ and $D^{\bullet}$. If we crush the 
boundary circles to points then we obtain  elements of  $\M(D^*)$, and there is a well defined
homomorphism $\T\rightarrow \M(D^*)$. 

\vspace{0.2cm}\noindent 
We studied in \cite{FK2} the asymptotic mapping class 
group $\M(D^*)$ denoted by $T^*$ there. Recall from \cite{FK2} that: 

\begin{proposition}
The group $T^*$ is generated by $\alpha$ and $\beta$. 
\end{proposition}

\noindent This implies that $\T\rightarrow T^*$ is an epimorphism. The relation between the asymptotic mapping class groups 
$\T$  and $T^*$ is made precise by the
following:

\begin{lemma}
We have an exact sequence 
\[ 0\to \Z^{\infty}\to \T \to T^*\to 1\]
where the central factor $\Z^{\infty}$ is the group of Dehn twists 
along attaching circles, normally generated by $t_0$. 
\end{lemma}
\begin{proof}
If an asymptotically rigid homeomorphism of $\su$ preserving $D^{\bullet}$ is isotopically trivial once the circles 
are crushed to points, then it is isotopic to a finite product of Dehn twists along those circles. Therefore, the kernel 
of $\T \to T^*$ is contained in the subgroup denoted $\Z^{\infty}$. Observe that $\alpha^4=t_0$ in $\T$, so that $t_{0}$ belongs to the
 kernel of $\T \to T^*$. Consequently, 
the kernel contains all the $\T$-conjugates of $t_{0}$, hence $\Z^{\infty}$. 
\end{proof}

\vspace{0.2cm}\noindent 
Thus if we understand the action of $T^*$ on  the isotopy classes of arcs 
embedded in $D^*$ then we can easily recover the action 
of $\T$ on  homotopy classes of loops of $\s^{\bullet}$, 
up to some twists along attaching circles. 

\subsubsection{The action of $T^*$ on the isotopy classes of arcs of $D^*$}

The planar model of $D^*$ is the punctured thick tree obtained from 
the binary tree by thickening in the plane and puncturing along midpoints 
of edges. The traces on $D^*$ of the loops coming from the pants decomposition of $\s$ are
 arcs transversal to the edges. Thus $D^*$ has a canonical 
decomposition into punctured hexagons. Each hexagon has three punctured 
sides coming from the arcs above, that we call {\it separating side arcs}. 
Moreover there are also three sides which are part of the boundary of 
$D^*$ that we will call {\it bounding side arcs}. Notice that hexagons correspond to 
vertices of the binary tree, while separating side arcs. Further $\beta$ is the rotation of order 3 supported on the 
hexagon $\overline{B}$ (image of the pants $B$) and $\alpha$ is the rotation of order  
4 that is supported on the union  of $\overline{B}$ with an adjacent hexagon.  
 
\begin{center}
\includegraphics{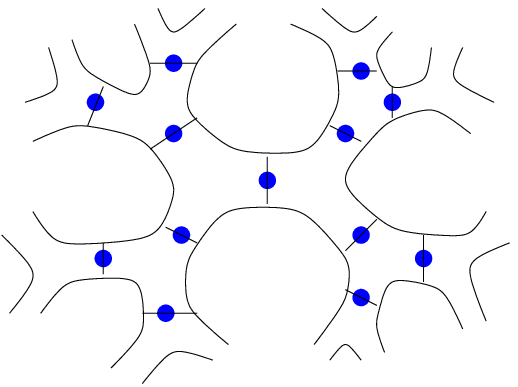}
\end{center}

\begin{lemma}\label{arctran}
Let $\gamma$ be an arc embedded in $D^*$ that joins two punctures. 
Then there exists some  element  of $T^*$ that sends $\gamma$ 
in a prescribed arc joining the punctures 0 and 1. 
\end{lemma}
\begin{proof}
Recall from \cite{FK2} that the infinite braid group associated to 
the punctures $B_{\infty}$ is contained in $T^*$. Further, there exists 
always a braid mapping class (supported in a compact subsurface of $D^*$) 
sending the arc $\gamma$ in the prescribed one. 
\end{proof}

\begin{lemma}\label{tranzsep}
The group $T^*$ acts transitively on the set of separating side arcs. 
\end{lemma}
\begin{proof}
The group $T^*$ contains $PSL_{2}(\Z)$, the group of orientation-preserving automorphisms
of the tree $\cal T$, generated by $\alpha^2$ and $\beta$. It acts transitively on the set of
edges of $\cal T$, hence on the set of separating sides of the hexagons of $D^*$.

\end{proof}

\noindent 
An arc joining a puncture belonging to a hexagon $H$ to a bounding side of $H$ is called 
{\it standard} if it is entirely contained in $H$.
\begin{lemma}
For any arc joining a puncture to a bounding side arc of a
hexagon, there exists 
some element of $T^*$ sending it into a standard arc joining the puncture 
$0$ to one of the bounding side of its hexagon. 
\end{lemma}
\begin{proof}
As above, $PSL(2,\Z)\subset T^*$ also acts transitively on the set of all bounding sides. 
Thus we can use an element of $T^*$ to send one end of our arc 
on a bounding side of the hexagon $\overline{B}$. Next, one composes by
 a braid element in $B_{\infty}$ that moves 
the other endpoint of the arc  
onto the puncture 0 and then makes the arc isotopic to a standard arc.  
\end{proof}

\noindent 
Let $t_{a_1}$ and $t_{b_1}$ denote the Dehn twists along a 
meridian $a_1$ and a longitude $b_1$ on the wrist.  

\begin{corollary}
The elements $\alpha,\beta, t, t_{a_1}, t_{b_1}, \pi$, a Dehn twist 
along one wrist connecting loop and a Dehn twist along a vertex connecting 
loop generate $\M$. 
\end{corollary}





\section{The rational homology of $\M$} 
 
\begin{theorem}\label{stab}      
The rational homology of $\M$ is isomorphic to the stable rational 
homology of the      
mapping class group: $H_*(\M,\Q)\cong  H_*(P\M,\Q)$.      
\end{theorem}      
      
\begin{proof}      
Recall first the theorem of stability, due to J. Harer (see \cite{Ha}): 
Let $R$ be a connected      
subsurface of genus $g_R$ of a connected compact surface $S$ with at least one boundary      
component. Then the map $H_n(P\M_R,\Z) \rightarrow H_n(P\M_S,\Z)$ induced by      
the natural morphism $P\M_R\rightarrow P\M_S$ is an isomorphism if $g_R\geq      
2n+1$.

\vspace{0.2cm}
\noindent      
The pure mapping class group $P\M$ is the inductive limit of the pure      
mapping class groups $P\M_R$, for all the compact subsurfaces $R\subset      
\su$. It follows that      
$H_n(P\M,\Q)=\displaystyle{\lim_{\stackrel{\rightarrow}{R}} H_n(P\M_R,\Q)}=$      
$H_n(P\M_R,\Q)$ for any compact subsurface $R\subset \su$ of genus $g_R \geq      
2n+1$. Therefore, the homology of $P\M$ is what is called the {\it stable homology of the mapping
class group}. By Mumford's conjecture proved in \cite{ma-we}, $H^*(P\M,\Q)$ is isomorphic
to $\Q[\kappa_{1},\ldots, \kappa_{i},\ldots]$, where $\kappa_{i}$,  the $i^{th}$ Miller-Morita-Mumford
class, has degree $2i$. Since $H^*(P\M,\Q)={\rm Hom}(H_{*}(P\M,\Q),\Q)$, each $H_{n}(P\M,\Q)$ is finite
dimensional over $\Q$.\\      
      
\noindent Write now the Lyndon-Hochschild-Serre spectral sequence in homology associated with  \[ 1      
\to P\M \to \M \to V\to 1\] The second term is      
$E^2_{p,q}=H_p(V,H_q(P\M,\Q))$. If we prove that $V$ acts trivially on the      
finite dimensional $\Q$-vector space $H_q(P\M,\Q)$, and invoke a theorem of      
K. Brown (\cite{br2}) saying that $V$ is rationally acyclic, then the only      
possibly non-trivial term of the spectral sequence is $E^2_{0,n}=      
H_n(P\M,\Q)$, and the proof is done.\\

\noindent Thus it remains to justify that  $V$ acts trivially on the homology groups $H_q(P\M,\Q)$, for any      
integer $q\geq 0$.  This results from the fact that $V$  is not linear, as we explain below. Indeed,   
if $\dim_{\Q}H_q(P\M,\Q)=N$, then ${\rm Aut}(H_q(P\M,\Q))\cong GL(N,\Q)$. So, let $\rho:V\rightarrow GL(N,\Q)$ be the representation      
resulting from the action of $V$ on $H_q(P\M,\Q)$. Since $V$ is a simple      
group, $\rho$ is either trivial or injective. Suppose it is injective, so that $V$  is
isomorphic to a finitely generated subgroup of $SL(N,\Q)$. Now each finitely generated
subgroup of  $SL(N,K)$ for any field $K$ is residually finite. But $V$ is not residually 
finite, since its unique normal subgroup of finite index is the trivial subgroup. Therefore, $\rho$
is trivial.

\end{proof}      
      
\begin{proposition}\label{c1} 
The free universal mapping class group $\M$ is perfect, and 
$H_2(\M,\Z)=\Z$. The generator of $H^2(\M,\Z)\cong \Z$ is called the first 
universal Chern class of $\M$, and is denoted $c_1(\M)$. 
 
\end{proposition}

\begin{proof} 
Recall that the pure mapping class group of a surface of type $(g,n)$ is 
perfect if $g\geq 3$. Consequently, $P\M$ is perfect. Since $V$ is perfect, 
$\M$ is perfect  as 
well.

\vspace{0.2cm}
\noindent
The above spectral sequence may be written with integral coefficients. One 
obtains $E^2_{2,0}=H_2(V,\Z)=0$ (see \cite{br2}), 
$E^2_{1,1}=H_1(V,H_1(P\M,\Z))=0$ since $P\M$ is perfect, and 
$E^2_{0,2}=H_0(V,H_2(P\M,\Z))$. By Harer's theorem (\cite{ha2}) and stability 
(\cite{Ha}), $H_2(P\M,\Z)\cong \Z$. The action of $V$ on  $\Z=H_2(P\M,\Z)$ 
must be trivial, since $V$ is simple, and it follows that 
$E^2_{0,2}=\Z$. Thus, the only non-trivial $E^{\infty}$ term is 
$E^{\infty}_{0,2}=E^2_{0,2}=\Z$, and this implies $H_2(\M,\Z)=\Z$. 
\end{proof}

\section{The symplectic representation in infinite genus} 

\subsection{Hilbert spaces and symplectic structure associated to ${\cal 
    S}_{\infty}$} 

There  is a natural intersection form 
$\omega:H_1({\cal S}_{\infty},\R)\times H_1({\cal   S}_{\infty},\R)\rightarrow \R$ on the homology of the infinite surface, but this is degenerate 
because it is obtained as a limit of intersection forms on 
surfaces with boundary. The ${\cal M}$-module $H_1({\cal S}_{\infty},\R)$ is the direct sum of two
submodules: $H_1({\cal S}_{\infty},\R)= H_1({\cal S}_{\infty},\R)_{s}\oplus  H_1({\cal S}_{\infty},\R)_{ns}$, where $H_1({\cal S}_{\infty},\R)_{s}$ is
generated by the homology classes of separating circles of ${\cal 
  S}_{\infty}$, while $H_1({\cal S}_{\infty},\R)_{ns}$ is generated by the homology classes of non-separating circles of ${\cal 
  S}_{\infty}.$  The kernel $\ker \omega$ of 
 $\omega$ is $H_1({\cal S}_{\infty},\R)_{s}$, and the restriction of $\omega$ to  $H_1({\cal S}_{\infty},\R)_{ns}$
 is a symplectic form.\\
 
\noindent 
For each wrist torus  occurring in the construction of ${\cal S}_{\infty}$ (see Definition \ref{surfinfgenus}), we consider
the meridian $a_k$ and the longitude $b_k$, with intersection number $\omega(a_k,b_{k})=1$. Note that   
both collections $\{a_k,\;k\in\N\}$ and 
$\{b_k,\; k\in\N\}$ are invariant by the mapping class group $\T$, since the generators $\alpha$
 and $\beta$ rigidly map  a wrist onto a wrist. Moreover, these collections are almost invariant by $\M$, meaning
 that for each 
$g\in {\cal M}$, $g(\{a_k,\;k\in\N\})$  (respectively 
 $g(\{b_k,\;k\in\N\})$) coincides up to isotopy with 
$\{a_k,\;k\in\N\}$ (respectively $\{b_k,\;k\in\N\}$) 
for all but finitely many elements.  \\

\begin{center}
\includegraphics{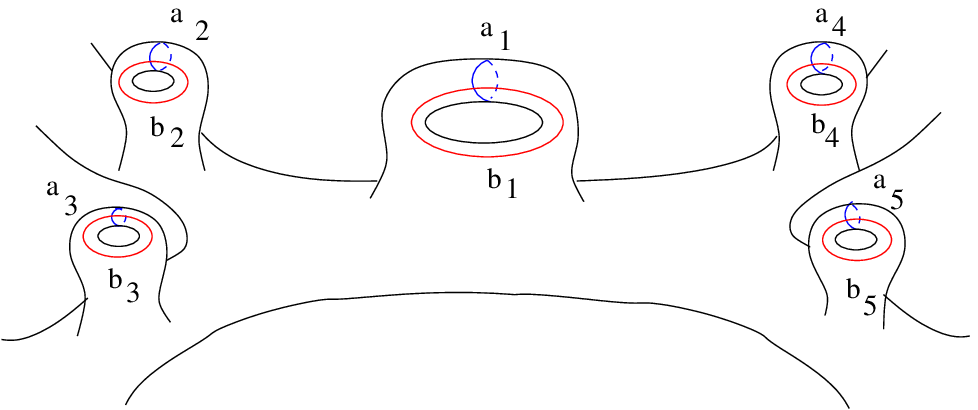}

\end{center}

\noindent The classes $\{a_k,b_k, k\in \N\}$ form a 
symplectic basis for  $H_1({\cal S}_{\infty},\R)_{ns}$. Each element of $\M$ acts on  $H_1({\cal S}_{\infty},\R)_{ns}$
by preserving the intersection form $\omega$. In particular, there is a representation
$$\rho: P\M \rightarrow {\rm Sp}(2\infty,\R)$$ 
where ${\rm Sp}(2\infty,\R)$ is the inductive limit of the symplectic groups ${\rm Sp}(2k,\R)$, with respect to
the natural inclusions ${\rm Sp}(2k,\R)\subset {\rm Sp}(2(k+1),\R)$. Note, though, that if $g\in\M$ is not in $P\M$, it is not represented into
${\rm Sp}(2\infty,\R)$, but into a larger symplectic group, that we are defining below.\\

\noindent One completes $H_1({\cal S}_{\infty},\R)_{ns}$ as
 a real Hilbert space for which this basis orthonormal. Let ${\cal H}_r$ be this Hilbert space, and $(.\,,.)$ denote its scalar product.\\ 
Let  $J$  be the almost-complex structure induced by $\omega$, i.e. 
the linear operator defined by  
$\omega(v,w)=(v,Jw)$ for all $v,w$ in ${\cal H}_r$. We have 
$J^2=-{\mathbf 1}$. 
Each  linear  operator  on ${\mathcal H}_r$ decomposes 
into a $J$-linear part $T_{1}$ and a $J$-antilinear part $T_{2}$, $T=T_1+T_2$, where 
$T_1=\frac{T-JTJ}{2}$ and $T_2=\frac{T+JTJ}{2}$. 

\vspace{0.2cm}\noindent 
Recall that (\cite{sh}) {\em the restricted symplectic group} ${\rm Sp_{res}}({\cal H}_r)$ 
is defined as the group of symplectic  (i.e. $\omega$ preserving)  
bounded invertible operators $T$ whose $J$-antilinear part $T_2$
is a Hilbert-Schmidt operator. An operator $T$ is called 
Hilbert-Schmidt if $||T||_{HS}^2:=\sum_i ||T(e_i)||^2$ is finite, 
where $(e_i)_{i\in\N}$ is an orthonormal Hilbert basis. 

\begin{theorem}\label{theo3}
The symplectic representation of the mapping class group $P\M$ extends to a representation $\hat{\rho}:\M\rightarrow  {\rm Sp_{res}}({\cal H}_r)$
of $\M$ into the restricted symplectic group.
\end{theorem}

\subsection{Proof of Theorem \ref{theo3}} 
 
Instead of a direct proof we will introduce the complexification of 
${\cal H}_r$ to be used  also in the next section. 
Let ${\cal H}={\cal  H}_r\otimes_{\R}\C$. 
Extend $\omega$ and $J$ by $\C$-linearity, and $(.\,,.)$ by  sesquilinearity, 
and denote by $\omega_{\C}, J_{\C}$ and $(.\,,.)_{\C}$ the  extensions. 
Thus, $({\cal H},(.\,.)_{\C})$ is  a complex Hilbert space. 
Let $ B$ be the indefinite hermitian form  
$ B(v,w)=\frac{1}{\sqrt{-1}}\omega_{\C}(v,\bar{w})$, 
for all $v,w$ in ${\cal H}$, where $\bar{w}$ is the complex-conjugate of $w$.

 \vspace{0.2cm}\noindent 
Let ${\rm Aut}({\cal H},\omega_{\C}, B)$ be the 
 group of bounded invertible operators of ${\cal H}$ which preserve 
 $\omega_{\C}$ and $B$. The morphism 
$\phi:{\rm Sp}({\cal H}_r)\longrightarrow {\rm Aut}({\cal H},\omega_{\C}, B)$, 
given by  $\phi(T)= T\otimes {\mathbf 1}_{\C}$ 
is an isomorphism (see \cite{ne}), since any  $T\in {\rm Aut}({\cal  H},\omega_{\C}, {\mathcal B})$ commutes with the complex conjugation and hence 
stabilizes ${\cal H}_r$.

\vspace{0.2cm}\noindent 
Since $J^2_{\C}=-id$, ${\cal H}={\cal H}_+\oplus {\cal H}_-$ where ${\cal 
  H}_{\pm}=\ker(J {\pm}\sqrt{-1}\cdot {\mathbf 1})$. Moreover, the direct sum is orthogonal. The 
complex conjugation interchanges ${\cal H}_+$ and ${\cal H}_-$. Let 
$(e_k)_{k\in\N}$ be an orthonormal basis of ${\cal H}_+$ and  
$(f_k=\overline{e_k})_{k\in\N}$ the conjugate basis of ${\cal H}_-$.  

\vspace{0.2cm}\noindent  
According to (\cite{ne}, 6.2) a symplectic operator $T$ belongs to ${\rm Sp_{res}}({\cal H}_r)$ if and only if the decomposition of $\phi(T)$ relative to 
the direct sum ${\cal H}_+\oplus {\cal H}_-$ in the basis $(e_k)_{k\in\N}\cup 
(f_k=\overline{e_k})_{k\in\N}$ reads   
$\left(\begin{array}{cc} 
\Phi & \Psi\\ 
\overline{\Psi}& \overline{\Phi} 
\end{array}\right)$, 
where 
\begin{enumerate} 
\item $^t{\overline\Phi}\Phi -{^t\Psi}{\overline\Psi}=1$ and 
 $^t{\overline\Phi}\Psi ={^t\Psi}{\overline\Phi}$, where 
$^tT$ denotes the adjoint of $T$ with respect to $(.\,,.)_{\C}$.  
\item $\Psi:{\mathcal H}_-\to {\mathcal H}_+$ is a Hilbert-Schmidt operator. 
\end{enumerate} 
We will apply this criterion for the action of ${\mathcal M}$.

\vspace{0.2cm}\noindent 
Set  
$$e_k=\frac{1}{\sqrt{2}}(a_k-\sqrt{-1}\cdot b_k), \,   
f_k=\frac{1}{\sqrt{2}}(a_k+\sqrt{-1}\cdot b_k),$$ 
Then $(e_k)_{k\in\N}$ is an orthonormal basis of ${\cal H}_+$ and 
$(f_k)_{k\in\N}$ is the conjugate orthonormal basis of ${\cal H}_-$.
Moreover, 
$\omega_{\C}(e_k,e_l)=\omega_{\C}(f_k,f_l)=0$, 
$\omega_{\C}(e_k,f_l)=\sqrt{-1}\delta_{kl}$, and 
${ B}(e_k,e_l)=-{ B}(f_k,f_l)=\delta_{kl}$, ${B}(e_k,f_l)=0$ for 
all $k,l$.

\vspace{0.2cm}\noindent 
Consider now  the action $\hat{\rho}(g)$ of 
$g\in {\cal M}$ on the ${\cal M}$-invariant subspace ${\cal H}_r$.  
We must check that $\Psi(\phi(\hat{\rho}(g)))$ is a Hilbert-Schmidt 
operator. In fact, it is a finite rank operator. 
Let $\Sigma_{h,n}$ be an admissible surface for 
$g$, that is $g$ is the mapping class of a homeomorphism  $G$ 
so that  $G: {\cal S}_{\infty} \setminus 
\Sigma_{h,n}\rightarrow  {\cal S}_{\infty} \setminus  \varphi(\Sigma_{h,n})$ 
is rigid. Any wrist torus $T_k=\Sigma_{1,1}$ of ${\cal S}_{\infty}\setminus \Sigma_{g,n}$ is rigidly mapped by $G$ onto another corresponding 
wrist torus $T_{\sigma(k)}$, for some infinite permutation $\sigma$. 
Therefore, for any such $T_k$, the associated matrices are such that
\[ \phi(\hat{\rho}(g))(e_k)=e_{\sigma(k)}, \,\,  \phi(\hat{\rho}(g))(f_k)=f_{\sigma(k)}\]
In particular, for all but finitely many $f_k$ (i.e.  excepting 
those corresponding to tori $T_k\subset \Sigma_{h,n}$) we have 
$\phi(\hat{\rho}(g))(f_k)\in {\cal H}_-$. Now 
$\Psi(\phi(\hat{\rho}(g))):{\cal H}_-\to {\cal H}_+$ corresponds 
to the components of $\phi(\hat{\rho}(g))(f_k)$ in ${\cal H}_+$. 
This means that $\Psi(\phi(\hat{\rho}(g)))$ has finite rank, and in 
particular, it is Hilbert-Schmidt. This proves that 
$\hat{\rho}(g)\in {\rm Sp_{res}}({\cal H}_r)$, as claimed.

\section{The universal first Chern class}
 
\subsection{The Pressley-Segal extension} 
 
Let 
${\cal H}$ be a polarized separable Hilbert space as above, that is, the orthogonal sum 
of two separable Hilbert spaces  ${\cal H}={\cal H}_+\oplus{\cal 
  H}_- $. 
The {\em restricted linear group} ${\rm GL_{res}}({\cal H})$ (see \cite{sh,pr-se}) 
is the Banach-Lie group of operators in ${\rm GL}({\cal H})$ 
whose block decomposition 
$A=\left( \begin{array}{cc} 
a & b\\ 
c & d 
\end{array} 
\right)$  
is such that $b$ and $c$ are Hilbert-Schmidt operators. 
Moreover, the invertibility of $A$ implies that $a$ is Fredholm in 
${\cal H}_+$, and has an index ${\rm ind}(a)\in \Z$. 
This gives a homomorphism 
${\rm ind}: {\rm GL_{res}}({\cal H})\rightarrow \Z$,  
that induces an isomorphism $\pi_0 
({\rm GL_{res}}({\cal H}))\cong \Z$. Denote by ${\rm GL^0_{res}}({\cal H})$ the connected 
component of the identity. Then ${\rm GL^0_{res}}({\cal H})$ is a perfect group 
(cf. \cite{co-ka}, \S5.4).

\begin{proposition} 
The restricted symplectic group ${\rm Sp_{res}}({\cal H}_r)$ embeds into the 
restricted linear group ${\rm GL^0_{res}}({\cal H})$. It is given the induced topology.  
\end{proposition} 
 
\begin{proof} 
Let $\left(\begin{array}{cc} 
\Phi & \Psi\\ 
\overline{\Psi}& \overline{\Phi} 
\end{array}\right)$ be in ${\rm Sp_{res}}({\cal H}_r)$. Since $\Psi$ is a 
Hilbert-Schmidt operator, $K={^t\Psi}{\overline\Psi}$ is 
trace-class, hence compact. Then $^t{\overline\Phi}\Phi= 1+K\geq 1$, hence 
$^t{\overline\Phi}\Phi\geq 1$ is injective, and the Fredholm alternative implies 
it is invertible. In particular, $\Phi$ itself is invertible, and has null index. 
\end{proof} 
 
\vspace{0.2cm}\noindent 
\noindent{\bf Pressley-Segal's extension of the restricted linear group.} Let ${\cal L}_1({\cal H}_+)$ denote the ideal of trace-class operators of 
${\cal H}_+$. It is a Banach algebra for the norm 
$||b||_1=Tr(\sqrt{b^*b})$, where $Tr$ 
denotes the trace form. We 
say that an invertible operator $q$ of ${\cal H}_+$ has a determinant if 
$q-id_ {{\cal H}_+}=Q$ is trace-class. Its determinant is the complex number 
$\det(q)=\sum_{i=0}^{+\infty} Tr(\wedge^i Q)$, where $\wedge^i Q$ is the 
operator of the Hilbert space $\wedge^i {\cal H}_+$ induced by $Q$ (cf. \cite{si}).

\vspace{0.2cm}\noindent  
Denote by ${\mathfrak T}$ the subgroup of $GL({\cal H}_+)$ consisting of operators 
which have a determinant, and by ${\mathfrak T}_1$ the kernel of the morphism 
${\rm det}: {\mathfrak T}\rightarrow \C^*$.

\vspace{0.2cm}\noindent 
Let ${\mathfrak E}$ be the subgroup of 
${\rm GL_{res}}({\cal H})\times {\rm GL}({\cal H}_+)$ 
consisting of  pairs $(A,q)$ such that  $a-q$ is trace-class. Then 
${\rm ind}(a)={\rm ind}(q+(a-q))={\rm ind}(q)=0$, so that $A$ belongs to 
${\rm GL^0_{res}}({\cal H})$. There is a short exact sequence 
$$\hspace{1cm} 1\rightarrow {\mathfrak T}\stackrel{i}{\longrightarrow} {\mathfrak E}\stackrel{p}{\longrightarrow} 
{\rm GL^0_{res}}({\cal H})\rightarrow 1$$ 
called the {\em Pressley-Segal extension}. 
Here $p(A,q)=A$, and $i(q)=({\mathbf 1}_{\cal H}, q)$. 
It induces the central extension  
  
$$1\rightarrow \frac{{\mathfrak T}}{{\mathfrak T}_1 }\cong \C^* \longrightarrow 
\frac{{\mathfrak E}}{{\mathfrak T}_1} \longrightarrow {\rm GL^0_{res}}({\cal H})\rightarrow 
1$$ 
The corresponding cohomology class in 
$H^2({\rm GL^0_{res}}({\cal H}),\C^*)$ is denoted by $PS$, and called the 
{\em Pressley-Segal class} of the restricted linear group.

\vspace{0.2cm}\noindent 
\noindent{\bf The Pressley-Segal class and the universal first Chern class.} 
For $G$ a Banach-Lie group, set $Ext(G,\C^*)$ for the set of equivalence 
classes of central extensions of $G$ by $\C^*$, 
which are locally trivial fibrations (see \cite{se}). 
Note that $Ext(G,\C^*)$ must not be confused with the 
group of continuous cohomology $H^2_{cont}(G,\C^*)$, since the latter only 
classifies the topologically split central extensions. One introduces two maps 
$H^2(G,\C^*)\stackrel{\delta}{\longleftarrow} 
Ext(G,\C^*)\stackrel{\tau}{\longrightarrow} H^2_{top}(G,\Z)$. The map $\delta$ 
associates to an extension $E$ in $Ext(G,\C^*)$ its cohomology class in the 
Eilenberg-McLane cohomology of $G$. The map $\tau$ is the composition of 
$Ext(G,\C^*) \longrightarrow [G,B\C^*=K(\Z,2)]$, which sends an extension to the 
homotopy class of its classifying map $G\rightarrow B\C^*$, with the 
isomorphism $[G,K(\Z,2)]\cong H^2_{top}(G,\Z)$.

\vspace{0.2cm}\noindent 
\noindent
Let us apply this formalism to $G={\rm GL^0_{res}}({\cal H})$ and the 
central Pressley-Segal extension, 
viewed as an element ${\cal PS}\in Ext({\rm GL^0_{res}}({\cal H},\C^*))$. 
Then $\delta({\cal PS})=PS$. The point 
is that ${\rm GL^0_{res}}({\cal H})$ is a homotopic model 
of the classifying space $BU$. In fact ${\mathfrak E}$ is contractible 
(see \cite{pr-se}, 6.6.2) and ${\mathfrak  T}$ is homotopically 
equivalent to $U$ (see \cite{pa}), hence the claim. 
It follows that the fibration ${\cal PS}$ corresponds to the universal 
first Chern class, that is, $\tau({\cal PS})=c_1(BU)\in H^2(BU,\Z)$.

\subsection{Cocycles on ${\rm GL^0_{res}}({\cal H})$, 
${\rm Sp_{res}}({\cal H}_r)$ and  $\M$} 
 \begin{lemma} 
The class $\iota^*(PS)$ in $H^2({\rm Sp_{res}}({\cal H}_r),\C^*)$ is represented by 
the cocycle 
$$C_1(g,g')=\det(\Phi(g)\Phi(g') \Phi(gg')^{-1})$$ 
\end{lemma}
\begin{proof}
Let ${\mathcal V}$ be the open subset  of 
${\rm GL^0_{res}}({\cal H})$ consisting of operators $A$ such that $a$ is 
invertible. It is known that the central Pressley-Segal extensions 
splits over ${\mathcal V}$, since it has the section $\sigma: 
{\mathcal V}\rightarrow {\mathfrak E},$ $A\mapsto(A,a)$. In particular, there is 
a local cocycle for $PS$ given by the formula (\cite{pr-se}, 6.6.4): 
$C(A,A')=\det(1+ aa'a''^{-1})$, for $A,A'\in {\mathcal V}$, 
where $a''$ is the first block of $A\cdot A'$. 
It suffices now to observe that ${\rm Sp_{res}}({\cal H}_r)$ 
embeds into ${\mathcal V}$. 
\end{proof}

In order to prove Proposition \ref{cocycle} below, we need the following 
result that  contrasts sharply with the finite dimensional case:
\begin{lemma} 
The restricted symplectic group ${\rm Sp_{res}}({\cal H}_r)$ is contractible. 
\end{lemma} 
\begin{proof} 
Denote by $Z$ the set of symmetric Hilbert-Schmidt operators ${\cal  H}_-\rightarrow {\cal  H}_+$ with norm $<1$. Clearly, $Z$ is a contractible 
subspace of the Banach space of Hilbert-Schmidt operators. 
The group ${\rm Sp_{res}}({\cal H}_r)$ acts transitively and 
continuously (see \cite{ne} p. 177) on $Z$ by means of 
$$g(S)=(\Phi(g) S+\Psi(g))(\overline{\Psi(g)} 
S+\overline{\Phi(g)})^{-1}\in Z, \,{\rm for }~ g\in {\rm Sp_{res}}({\cal H}_r), 
S\in Z$$ 
The stabilizer of $S=0$ is the group of matrices $\left(\begin{array}{cc} 
\Phi & 0\\ 
0 & \overline{\Phi} 
\end{array}\right)$ 
such that $\Phi$ is unitary in ${\cal H}_+$. Thus, it is isomorphic to ${\cal 
  U}({\cal H}_+)$. By a result of Kuiper (\cite{ku}), ${\cal 
  U}({\cal H}_+)$ is contractible. The claim is now a consequence 
of the contractibility of  ${\rm Sp_{res}}({\cal H}_r )/{\cal 
  U}({\cal H}_+)\cong Z$. 
\end{proof}     
 
\begin{proposition}\label{cocycle} 
For each integer $n\in\Z$, there is a well-defined continuous cocycle $C_n$ defined on 
${\rm Sp_{res}}({\cal H}_r)$, with values in $\C^*$, such that 
$$C_n(g,g')=\det\left( (\Phi(g)\Phi(g')\Phi(gg')^{-1})^{\frac{1}{n}}\right)$$ 
Moreover, $\frac{C_n}{|C_n|}$ may be lifted to a real cocycle 
$\widehat{\varsigma_n}:{\rm Sp_{res}}({\cal 
  H}_r)\times {\rm Sp_{res}}({\cal  H}_r) \longrightarrow \R$ such that 
\[\frac{C_n(g,g')}{|C_n(g,g')|}=e^{2i\pi \widehat{\varsigma_n}(g,g')}, \,\, 
{\rm for ~ all~} g,g'\in {\rm Sp_{res}}({\cal  H}_r)\]  
The restriction $\varsigma_{1}$ of $\widehat{\varsigma_1}$ to  ${\rm Sp}(2\infty,\R)$ defines an integral cohomology
class $[\varsigma_{1}]\in H^2({\rm Sp}(2\infty,\R), \Z)$.

\end{proposition}

\begin{proof} 
In fact, $\Phi(g)^{-1} \Phi(g')^{-1}\Phi(gg')= 
1+(\Phi(g')^{-1}\Phi(g)^{-1}\Psi(g)(\overline{\Psi(g')})$. But, according to 
(\cite{ne}, p. 168) we have 
\[ ||\Phi^{-1}(g)\Psi(g)||<1 \;\;{\rm and }\;\; ||\overline{\Psi(g')}{\Phi(g')}^{-1}||<1 \]
Thus, there is a non-ambiguous definition of 
$(\Phi^{-1}(g) \Phi(gg')\Phi^{-1}(g'))^{\frac{1}{n}}$ given by an absolutely 
convergent series.\\

\noindent  
The existence of $\varsigma_n$ is now an immediate consequence of the 
preceding lemma. \\

\noindent The map $\ell: g\in {\rm Sp}(2\infty,\R)\mapsto \ell(g)=\frac{det(\Phi(g))}{\mid det(\Phi(g))\mid}$ 
is well-defined,  so that the cocycle 
$$(g,g')\in {\rm Sp}(2\infty,\R)\times {\rm Sp}(2\infty,\R)\mapsto e^{2i\pi \varsigma_{1}(g,g')}$$ is the
coboundary of $\ell$. This proves that the cohomology class of $\varsigma_{1}$ restricted to ${\rm Sp}(2\infty,\R)$
is integral.
\end{proof} 
 
\begin{remark}\label{dupont}

\begin{enumerate} 

\item  The restrictions of the real cocycles 
$\varsigma_n$ on the finite dimensional Lie 
group ${\rm Sp}(2g,\R)$ are those constructed by Dupont-Guichardet-Wigner (see 
\cite{gu-wi}). In fact, the authors of \cite{gu-wi} proved that the cohomology class of the restriction of  $\varsigma_1$ to ${\rm Sp}(2g,\R)$ 
is integral, and is the image in $H^2({\rm Sp}(2g,\R),\R)$ of 
the generator of $H_{\rm bor}^2({\rm Sp}(2g,\R),\Z)=\Z$, 
the second group of borelian cohomology of ${\rm Sp}(2g,\R)$. 
They prove also that it is the image of the first Chern class $c_{1}(BU(g,\C))$ by the composition of maps
$$H^2(BU(g,\C),\Z)\approx H^2(B {\rm Sp}(2g,\R),\Z)\rightarrow H^2(B {\rm Sp}(2g,\R)^{\delta},\Z)\approx H^2({\rm Sp}(2g,\R),\Z)\rightarrow H^2({\rm Sp}(2g,\R),\R),$$
where $B{\rm Sp}(2g,\R)^{\delta}$ is the classifying space of ${\rm Sp}(2g,\R)$ as a discrete group.\\

\item The remark above implies that the map
  $$H^*(BU,\Z)\approx H^*( B{\rm Sp}(2\infty,\R),\Z  )\rightarrow H^*({\rm Sp}(2\infty,\R),\Z  ).$$
  sends the first universal Chern class $c_{1}(BU)$ onto $[\varsigma_{1}]\in H^2({\rm Sp}(2\infty,\R), \Z).$
  Further, the symplectic representation $\rho:P\M\rightarrow {\rm Sp}(2\infty,\R)$ maps $[\varsigma_{1}]$ onto
  the generator $c_{1}(P\M)$ of $H^2(P\M,\Z)$.

 \item According to (\cite{ne}, Theorem 6.2.3), the Berezin-Segal-Shale-Weil cocycle is the complex conjugate of the cocycle 
$C_{-\frac{1}{2}}$.
   
\end{enumerate}
\end{remark} 
 
\begin{theorem}\label{chernimp} 
Let $[\hat{\varsigma_1}]\in H^2({\rm Sp_{res}}({\cal H}_r),\R)$ be the cohomology class of 
$\hat{\varsigma_1}$. The pull-back of $[\hat{\varsigma_1}]$ in $H^2({\cal M},\R)$ by the 
representation $\hat{\rho}$ of Theorem \ref{metab} is integral, and is the natural 
image of the generator $c_{1}(\M)$ of $H^2({\cal M},\Z)$ in $H^2({\cal M},\R)$. 
\end{theorem}

\begin{proof} 
   Let $\iota:P\M\rightarrow \M$  and $j:{\rm Sp}(2\infty,\R)\rightarrow {\rm Sp_{res}}({\cal H}_r)$   be
   the natural embeddings. Plainly, $j\circ \rho=\hat{\rho}\circ \iota$. Since $j^*: H^2({\rm Sp_{res}}({\cal H}_r),\R)\rightarrow H^2({\rm Sp}(2\infty,\R),\R)$ maps 
  $[\hat{\varsigma_1}]$ onto $[{\varsigma_1}]$, one has $\iota^*(\hat{\rho}^*[\hat{\varsigma_1}])=\rho^*[\varsigma_1]$. 
  Let us denote by $\bar{c_{1}}(P\M)$ (respectively $\bar{c_{1}}(\M)$) the image of ${c_{1}}(P\M)$ (respectively ${c_{1}}(\M)$) in $H^2(P\M,\R)$ 
(respectively $H^2(\M,\R)$).
   According to Remark \ref{dupont}, 2.,  $\rho^*[\varsigma_1]=\bar{c_{1}}(P\M)$. By Proposition \ref{c1}, $\iota^*(\bar{c_{1}}(P\M))= \bar{c_{1}}(\M)$,
   hence $\hat{\rho}^*[\hat{\varsigma_1}]= \bar{c_{1}}(\M)$.
\end{proof}

 {\small     
      
\bibliographystyle{plain}      
      
}      
\end{document}